\documentclass[reqno,a4paper,final]{amsart}
\usepackage{amssymb,amstext,amsthm,eucal}
\usepackage{bbold}
\usepackage{array}
\usepackage[latin1]{inputenc}
\newcommand{\SO}{\mathrm{SO}}
\newcommand{\CO}{\mathrm{CO}}

\newcommand{\Spin}{\mathrm{Spin}}

\newcommand{\SU}{\mathrm{SU}}

\newcommand{\RR}{\mathbb{R}}

\newcommand{\ZZ}{\mathbb{Z}}

\newtheorem{THEO}{\bf Theorem}
\newtheorem{DEF}{\bf Definition}
\newtheorem{PR}{\bf Proposition}

\newcommand{\MUNCH}[1]{\relax}
\allowdisplaybreaks[1]
\begin{document}
\begin{sloppypar}
\title{Conformal holonomy of bi-invariant metrics}
\author{Felipe Leitner}
\address{Institut f{\"u}r Mathematik, Universit{\"a}t Leipzig, Germany}
\email{leitner@mathematik.uni-leipzig.de}
\thanks{ }
\date{\today}

\begin{abstract} 
We discuss in this paper the conformal geometry of bi-invariant
metrics on compact semisimple Lie groups. For this purpose
we develop a conformal Cartan calculus adapted to this problem.
In particular, we derive an explicit formula for the holonomy
algebra of the normal conformal Cartan connection of a bi-invariant
metric. As an example, we apply this calculus to the group $\SO(4)$. Its 
conformal
holonomy group is calculated to be $\SO(7)$.  
\end{abstract}
\maketitle
\tableofcontents

\section{Introduction}
\label{ab1}

In Riemannian geometry the concept of holonomy is a well studied
problem. The holonomy groups arise from affine connections on
Riemannian manifolds. Most famous is the holonomy theory of 
the Levi-Civita connection, which is the canonical
connection to a Riemannian metric. A list of possible 
holonomy groups in this case was first established in \cite{Ber55}.

In conformal geometry, there exists no canonical affine connection. Any
choice of a torsion-free affine connection preserving the structure
group (Weyl connection) is an additional datum to the conformal 
structure on a manifold. However, taking into consideration that
conformal geometry should be understood as a 'second order' structure,
it turns out and is a well-known fact that there exists a canonical
Cartan connection on the prolongated principal fibre bundle of
second order, which has a parabolic of the M{\"o}bius group
as structure group. This canonical connection solves Cartan's equivalence
problem for conformal geometry, which says that  two conformal 
structures
are equivalent if and only if the canonical connections coincide.
The canonical connection on a conformal space takes in general values in 
the Lie algebra $\frak{so}(1,n+1)$ of the M{\"o}bius group $\SO(1,n+1)$. 
As in the case of
Riemannian geometry,
the knowledge of the holonomy algebra and group defined to this canonical
conformal Cartan connection is a significant invariant of a space
with conformal structure. In particular, it describes invariant 
structures and solutions of invariant
differential operators on a conformal space.

The aim of this paper is to develop and apply a conformal Cartan
calculus to a very simple situation, namely the conformal geometry
of bi-invariant metrics on compact semisimple Lie groups, in order
to calculate the conformal holonomy explicitly.
We proceed as follows. In the next two paragraphs, we recall briefly the 
basic concept of 
conformal 
Cartan geometry, in general, and the well-known
notion of bi-invariant metrics. In paragraph
\ref{ab4}, we develop the conformal Cartan calculus on bi-invariant
metrics. In particular, we will describe the canonical connection
and its curvature on a semisimple group $N$ by certain maps 
$\gamma_{nor}$
and $\kappa$, which live on 
the flat 
model 
$(\frak{so}(1,n+1),\frak{p})$. Then we discuss properties of these maps
and derive an explicit formula for the holonomy algebra (paragraph 
\ref{ab5}). Finally, we make explicit calculations for the bi-invariant
metric on $\SO(4)$, in particular, we derive the holonomy group.

It is the project of a forthcoming paper (cf. \cite{Lei04b}) to 
investigate
the conformal holonomy in more general situations e.g. like   
invariant metrics on (reductive) homogeneous Riemannian spaces. 
Even more generally, it is the idea to find explicit formulas 
and calculations for the holonomy of homogeneous parabolic geometries 
admitting a
canonical Cartan connection (remember that conformal geometry is an 
example for a 
$|1|$-graded
parabolic geometry). With 
homogeneous space we do not mean here
a space with flat canonical connection, but a space with a parabolic 
geometry whose automorphism
group acts transitively on the base space. Beside conformal geometry, a 
further classical example for this approach would be CR-geometry, which 
is a $|2|$-graded parabolic geometry.

\section{Conformal Cartan geometry}
\label{ab2}

We describe here briefly the conformal structure of a smooth
manifold uniquely determined by a normal Cartan connection on a second 
order principal
fibre bundle with parabolic structure group (cf. \cite{Kob72},
\cite{CSS97}). We 
start with 
recalling the 
flat homogeneous model of conformal geometry.

Let $G=\SO(1,n+1)$ be the Lorentzian group, which acts on the
$(n+2)$-dimensional Minkowski space $\RR^{1,n+1}$ equipped with the
scalar product
\[\langle x,x\rangle_{1,n+1}= 2x_0x_{n+1}+\sum_{i=1}^{n}x_i^2\ .\]
Its Lie algebra $\frak{g}=\frak{so}(1,n+1)$ is $|1|$-graded by
\[\frak{g}=\frak{m}_{-1}\oplus\frak{co}(n)\oplus \frak{m}_1\ ,\]
where $\frak{m}_{-1}=\RR^n$ and $\frak{m}_1=\RR^{n*}$ are dual 
vector spaces
via the Killing form on $\frak{g}$. The 
$0$\,-part
$\frak{p}_0:=\frak{co}(n)$ of this grading decomposes to the
semisimple part $\frak{so}(n)$ and the center $\RR$, which is responsible
for the conformal weight of a representation of $\frak{g}$. In matrix 
form, we have 
\[\left(\begin{array}{ccc} 0& 0& 0\\m & 0& 0\\ 0& -{}^tm& 0
\end{array}\right)\in\frak{m}_{-1}\ ,\quad\ \ 
\left(\begin{array}{ccc} -a& 0& 0\\0 & A& 0\\ 0& 0& a
\end{array}\right)\in\frak{p}_0\ ,\quad\ \
\left(\begin{array}{ccc} 0& l& 0\\0 & 0& -{}^tl\\ 0& 0& 0
\end{array}\right)\in\frak{m}_{1}\ . \] 
The commutators are then given by
\[ \begin{array}{ll}
{[\ ,\ ]}:\ \frak{p}_0\times\frak{p}_0\to \frak{p}_0\ ,\qquad & 
[(A,a),(A',a')]=(AA'-A'A,0)\\
{[\ ,\ ]}:\ \frak{p}_0\times \frak{m}_{-1}\to \frak{m}_{-1}\ ,\qquad &
[(A,a),m]=Am+am\\
{[\ ,\ ]}:\ \frak{m}_1\times \frak{p}_{0}\to \frak{m}_{1}\ ,\qquad &
[l,(A,a)]=lA+al\\
{[\ ,\ ]}:\ \frak{m}_{-1}\times \frak{m}_{1}\to \frak{p}_{0}\ ,\qquad &
[m,l]=(ml-{}^t(ml),lm)\ ,\end{array}\qquad \]
where $(A,a)$, $(A',a')\in \frak{so}\oplus\RR,\ m\in\RR^n,\ l\in
\RR^{n*}$.
The subalgebra
\[\frak{p}:=\frak{co}(n)\oplus\frak{m}_1\]
is parabolic, i.e., it contains a maximal solvable subalgebra of
$\frak{so}(1,n+1)$.

Let $P$ be the closed subgroup of $G$ consisting of the 
matrices
\[ \left\{\left(\begin{array}{ccc} a^{-1}& v& b\\0 & A& r\\ 0& 0& a
\end{array}\right)\left|\begin{array}{c} A\in SO(n),\ 
a\in\RR\smallsetminus 0, 
v\in\RR^{n*},\\[1mm]
r=-aA{}^tv,\ b=-\frac{a}{2}\langle {}^tv,{}^tv\rangle_n\end{array}
\right.\right\}\ .\]
The group $P$ is the stabilizer of the point
$o=[1:0:\cdots:0]$ 
in the projectivation $P_{n+1}(\RR)$ of the Minkowski space $\RR^{1,n+1}$
and the Lie algebra of $P$ is the parabolic $\frak{p}$. The homogeneous 
space
$G/P$ is the $n$-dimensional M{\"o}bius sphere $S^n$ with standard 
conformally flat structure, which is induced here from the 
scalar product on $\RR^{1,n+1}$, and $G$ acts by conformal 
automorphisms on $S^n$. Furthermore,
let $P_1$ be the closed subgroup of $P$ given by the matrices
\[ \left\{\left(\begin{array}{ccc} 1& v& b\\0 & I_n& -{}^tv\\ 0& 0& 1
\end{array}\right)\left|\
v\in\RR^{n*},\ b=-\frac{1}{2}\langle 
{}^tv,{}^tv\rangle\right.\right\}\ 
.\]
This is the vector group $\exp \frak{m}_1$ with Lie 
algebra 
$\frak{m}_1$ and it is the kernel of the linear isotropy representation
of $P$ acting on $T_oS^n$ (the point $o$ can be thought of as the point 
at 
infinity when conformally compactifying the Euclidean space 
$\RR^n$).  
Since $P_1$ is normal, the quotient $P/P_1$ is a group itself, which is
isomorphic
to $\CO(n)$. Any 
element $p\in P$ can 
be represented in a unique way 
by a product $p_0\cdot \exp l$, where $p_0\in \CO(n)$ and $l\in 
\frak{m}_1$.  
   
The Maurer-Cartan form 
\[\omega_G: TG \to \frak{g}\]
is a Cartan connection on $S^n$ with parabolic structure  group $P$.
The Maurer-Cartan equation 
\[ d\omega_G=-\frac{1}{2}[\omega_G,\omega_G]\  \]
shows that the curvature of  the Cartan connection $\omega_G$ vanishes,
i.e., $\omega_G$ is a flat connection.

We apply now the flat model $(G,P)$ to (conformally) curved spaces. 
Let $(M^n,c)$ be a smooth manifold of dimension $n\geq 3$ with
conformal structure $c$. The conformal structure $c$ 
is usually given  
as an equivalence class of metrics $[g]$, which differ from each other 
only by multiplication 
with a positive smooth function on $M$ (scaling function). The conformal 
structure $c$ defines
a first order reduction $CO(M)$ of the general linear frame bundle $GL(M)$ 
to the structure group $\CO(n)$ in the usual manner. On the other side, 
the conformal 
structure also
induces a reduction $P(M)$ of the general linear frame bundle 
$GL^2(M)$ 
of second order to the structure group $P$. It is \[P(M)/P_1
\cong CO(M)\ .\] Then the 
$\CO(n)$-invariant
lifts from $CO(M)$ to $P(M)$ are in $1$-to-$1$ correspondence with
Weyl connections on $(M,c)$, i.e., affine connections 
with structure group $\CO(n)$ and without torsion. 

The principal $P$-fibre bundle $P(M)$ admits a canonical Cartan connection 
$\omega_{nor}$, which, in particular, describes the embedding of $P(M)$ 
into 
$GL^2(M)$, 
and therefore determines the conformal structure $c$ on $M$. The canonical
connection is made unique by imposing a normalisation condition.   
To say it in detail, the Cartan connection 
\[\omega_{nor}:TP(M)\to \frak{g}\] has the properties: 
\begin{enumerate}\item $\omega_{nor}(p): T_pP(M)\to \frak{g}$ is 
an 
isomorphism for all 
$p\in P(M)$, which is along the fibre the tautological map, and
\item $R_p^*\omega_{nor}=\mbox{Ad}(p^{-1})\omega_{nor}$, where 
$R_p$ denotes right translation by $p$ on $P(M)$.
\end{enumerate}
The curvature of $\omega_{nor}$ is defined by
\[K:=d\omega_{nor}+\frac{1}{2}[\omega_{nor},\omega_{nor}]\ ,\]
where 
\[[\omega_{nor},\omega_{nor}](X,Y):=[\omega_{nor}(X),\omega_{nor}(Y)]
-[\omega_{nor}(Y),\omega_{nor}(X)]\]
for all $X,Y\in TP(M)$. It holds $R_p^*K=\mbox{Ad}(p^{-1})K$ and  
$K$ is in vertical direction 
trivial, i.e., inserting a vertical vector on $P(M)$ into $K$ produces
zero. 
Alternatively, the curvature is uniquely determined by
the function \[\kappa=\kappa_{-1}\oplus\kappa_0\oplus\kappa_1: P(M) \to
\frak{m}_{-1}^*\otimes\frak{m}_{-1}^*\otimes\frak{g}\ ,\] 
which is defined by
\[\kappa_p(a,b)=K_p(\omega_{nor}^{-1}(a),\omega_{nor}^{-1}(b))\]
for all $p\in P(M)$ and  $a,b\in \frak{m}_{-1}$. The normalisation 
condition,
which makes the connection $\omega_{nor}$ unique, is expressed by
the curvature properties 
\begin{enumerate}
\item $\kappa_{-1}=0$, i.e., the connection $\omega_{nor}$ is 
torsion-free,
and 
\item the trace-free condition 
$tr\kappa_{0}=\sum_{i=1}^n\kappa_0(e_i,a)(b)(e_i^*)=0$ for all
$a,b\in \frak{m}_{-1}$.
\end{enumerate}  Thereby, the $e_i$'s denote the standard basis
of $\frak{m}_{-1}\cong\RR^n$ and the $e_i^*$'s are dual in 
$\frak{m}_{1}$.

The canonical normal Cartan connection and its curvature can be 
described with 
respect
to a metric $g$ in the conformal class $c$ as follows.
The metric $g$ determines a Weyl connection $\sigma_g$, i.e., an 
invariant 
lift 
from 
$CO(M)$ to $P(M)$. Then the pull-back $\sigma_g^*\omega_{-1}$ of the 
$(-1)$-part of 
$\omega_{nor}$  
gives the soldering form on $CO(M)$, which in turn identifies (with 
respect to a base frame) the tangent space
of $M$ at a point with $\frak{m}_{-1}$. The 
$0$-part $\sigma_g^*\omega_{0}$
is the Levi-Civita connection form to $g$ on $CO(M)$ and 
$-\sigma_g^*\omega_{1}$ is tensorial and projects to $(M,g)$ as the
Schouten tensor ('rho'-tensor), which is   
defined by
\[L(X)=\frac{1}{n-2}\big(\frac{scal}{2(n-1)}X-Ric(X)\big)\ ,\]
where $Ric$ denotes the Ricci tensor and $scal$ the scalar curvature
to $g$ on $M$. The (harmonic) $0$-part $K_0$ of the curvature 
of $\omega_{nor}$ corresponds to 
the
Weyl tensor $W$ on $(M,g)$, which is the trace-free part of the Riemannian
curvature tensor $R$, and can be expressed by
\[W:= R-g*L\ ,\] 
where $g*L$ denotes a Kulkarni-Nomizu product. Finally, the negative 
$-K_1$ of 
the $1$-part of the curvature projects to the Cotton-York tensor $C$ on 
$(M,g)$, which is 
given in terms of $g$ by
\[C(X,Y):=\nabla^g_XL(Y)-\nabla^g_YL(X)\ .\]

\section{Bi-invariant metrics }                    
\label{ab3}                                        

We recall here the notion of bi-invariant metrics on
compact semisimple Lie groups. They will be the matter 
of our investigation when
we apply in the next paragraph the Cartan calculus to 
their conformal classes.

Let $N$ be a connected and compact semisimple Lie group of dimension 
$n$
and let $\frak{n}$ denote its Lie algebra. The Killing form
\[B(X,Y):=\mbox{tr}\, \mbox{ad}_X\mbox{ad}_Y\] is 
$\mbox{Ad}(N)$-invariant
and negative definite on $\frak{n}$.
In particular,
\[B(X,[Y,Z])=B([X,Y],Z)\qquad \mbox{for\ all\ } X,Y,Z\in\frak{n}\ .\]
The negative $-B$ of the Killing form defines through left translation
with the group multiplication an invariant metric $g_{\frak{n}}$
on $N$. In fact, the metric $g_{\frak{n}}$ is not only left-invariant,
but also right-invariant under the group multiplication and the metric
$g_{\frak{n}}$ is called a bi-invariant metric on $N$ (cf. 
\cite{O'N83}).
In the following, we will often identify 
left-invariant 
vector fields on $N$ with their generators in $\frak{n}$,
\[\tilde{X}(n):=\frac{d}{dt}|_{t=0}n\cdot\exp tX\quad\mapsto\quad 
X\in\frak{n}\ .\]
For the Levi-Civita connection $\nabla$ of $g_{\frak{n}}$, it holds
\[\nabla_XY=\frac{1}{2}[X,Y]\qquad\mbox{for\ all\ }\  X,Y\in \frak{n}\ 
.\]
The Riemannian curvature tensor is then given by
\[R_{XY}Z=-\frac{1}{4}[[X,Y],Z]\qquad\mbox{for\ all\ }\ X,Y,Z\in 
\frak{n}\]
and the sectional curvature of a plane spanned by orthonormal elements
$X,Y$ in $\frak{n}$
is 
\[S(X,Y):= -\frac{1}{4}B([X,Y],[X,Y])\ .\]
For the Ricci tensor we find
\[Ric|_\frak{n}=-\frac{1}{4}B\ ,\]
i.e., $Ric=\frac{1}{4}g_{\frak{n}}$ and $g_{\frak{n}}$ is an Einstein 
metric on $N$ with positive scalar curvature $scal=\frac{n}{4}$. The 
Schouten tensor is $L=-\frac{1}{8(n-1)}g_{\frak{n}}$ and  the 
Cotton-York
tensor $C$ vanishes identically, $C\equiv 0$. For the Weyl tensor, we 
obtain
\[W=R+\frac{1}{8(n-1)}g_{\frak{n}}*g_{\frak{n}}\ .\]

\section{Conformal Cartan geometry and bi-invariant metrics} 
\label{ab4}                                                 
At the end of the last paragraph we calculated already the content 
of the 
conformal
curvature, the Weyl tensor and the Cotton-York tensor, in terms of
the bi-invariant metric on a compact semisimple group. 
However, we want to establish here a
conformal Cartan calculus for bi-invariant metrics.
This will help us to get a better understanding of the normal
connection and its curvature. Though, our approach will always
use a convenient trivialisation, which represents the bi-invariant
metric in the conformal class.

Let $N$ be a connected and compact semisimple Lie group of dimension
$n$ with
Lie algebra $\frak{n}$ and bi-invariant metric $g_{\frak{n}}$.
Then there is the conformal structure $c_{\frak{n}}=[g_{\frak{n}}]$ 
defined on
the group $N$. Let 
\[\theta:(\frak{n},-B)\ \cong\ (\frak{m}_{-1},\langle\ ,\ \rangle_n)\]
be an isometry and $(e_1,\ldots,e_n)$ the standard basis on 
$\frak{m}_{-1}\cong\RR^n$. 
The map $\theta$ transfers the Lie bracket $[\ ,\ ]_{\frak{n}}$ to 
$\frak{m}_{-1}$ through the expression
\[\rho_{\frak{n},\theta}(a,b):=\theta 
[\theta^{-1}(a),\theta^{-1}(b)]_{\frak{n}}\ ,\]
where $a,b\in \frak{m}_{-1}$. (Likewise, we will also use 
the notation $[\ ,\ ]_{\frak{n}}$ on $\frak{m}_{-1}$).
Moreover, the map
$\theta$ induces the orthonormal frame 
\[ \left\{E_i:=\theta^{-1}(e_i)|\quad i=1,\ldots,n\right\}\]  
on $\frak{n}\cong T_eN$. The corresponding left-invariant frame field 
gives rise to a global trivialisation of the first order 
conformal frame bundle $CO(N)$ on $(N,c)$,
\[\begin{array}{ccc}CO(N)&\quad\cong\qquad& N\times\CO(n)\\
 \left\{\tilde{E}_i(s)|\ i=1,\ldots, n\right\}&\quad\mapsto\qquad& (s,e)
\end{array}\]
(in fact, it is a trivialisation of the orthonormal frame bundle
of $(N,g_{\frak{n}})$).
We denote by $P(N)$ the second order conformal frame bundle with
parabolic structure group $P$, which is a subbundle of $GL^2(N)$.
The bi-invariant metric $g_{\frak{n}}$ induces an invariant lift
\[ \sigma_{g_\frak{n}}:CO(N)\to P(N)\] and with the help of 
the $\tilde{E}_i$'s also a trivialisation 
\[\iota_{\theta,g_\frak{n}}:P(N)\quad\cong\quad N\times P\]
of the parabolic frame bundle. The left translation $L_s:N\to N$ 
preserves for all $s\in N$ the bi-invariant metric, hence the conformal
structure $c_{\frak{n}}$. The translation $L_s$ on $N$ induces in 
a natural
way transformations (also denoted by $L_s$) of the conformal frame 
bundles
$CO(N)$ and $P(N)$. As we have chosen here the trivialisation 
$\iota_{\theta,g_{\frak{n}}}$, it holds 
\[\begin{array}{cccl} 
\iota_{\theta,g_{\frak{n}}}\circ L_s\circ\iota_{\theta,
g_{\frak{n}}}^{-1}:& 
N\times P&\to & N\times P\ .\\
& (n,p)&\to& (sn,p)\end{array}\]

There exists a canonical Cartan connection on $P(N)$ denoted by
\[\omega_{nor}:TP(N)\to \frak{g}\ .\] This connection is determined by  
its curvature properties, namely the curvature  $K$ is torsion-free and 
satisfies 
the 
trace-free condition
on the $0$-part $K_0$ (cf. paragraph \ref{ab2}). By definition, the 
connection $\omega_{nor}$ is right-invariant along the fibres, i.e., 
\[R_p^*\omega_{nor}=\mbox{Ad}(p^{-1})\omega_{nor}\quad\ \mbox{for\ 
all\ }p\in P\ ,\]
and it is also left-invariant under the group multiplication on
$N$, i.e., 
\[L_s^*\omega_{nor}=\omega_{nor}\ .\]
The latter fact is, because, 
the normal Cartan connection is uniquely determined by
the conformal structure and is therefore invariant under the conformal
automorphism group. These two facts, the right-invariance by
multiplication with $P$ and the left-invariance with respect to $N$, imply 
that
$\omega_{nor}$ on $P(N)$ is uniquely determined by the 
(vector space) isomorphism induced from $\omega_{nor}$ at a single point 
$x_o$ of $P(N)$. Using the 
trivialisation $\iota_{\theta,g_{\frak{n}}}$ from above, 
we can choose this point 
to 
$x_o=(e,e)$ and 
then the
isomorphism  
\[\omega_{nor}(e,e): \frak{n}\times\frak{p}\to 
\frak{g}\]
determines the canonical Cartan connection.
Since the $(-1)$-part of $\omega_{nor}$ corresponds to the soldering form,
it holds $\omega_{-1}(e,e)(E_i)=e_i$. This leads us to the 
definition
of the map $\gamma_{nor}$ (which depends on the choice of 
$\theta$) through
\[\begin{array}{cccl}\gamma_{nor}\ :& \frak{m}_{-1}&\to&\frak{p}\ \ 
,\\
& a &\mapsto&
\pi_{\frak{p}}\circ\omega_{nor}(e,e)\circ\pi_{\frak{n}}\circ
\omega_{nor}^{-1}(e,e)(a)
\end{array}\]
where $\pi_{\frak{p}}$ and $\pi_{\frak{n}}$ denote the obvious projections
(the latter with respect to our trivialisation). This map decomposes
to \[\gamma_{nor}=\gamma_0+\gamma_1\] and
it still contains the whole information of
the canonical Cartan connection $\omega_{nor}$ on $N$, since it holds
the relation
\[\omega_{nor}(E_i)= e_i+\gamma_{nor}(e_i)\ .\]
In fact, $\omega_{nor}$ is recovered from $\gamma_{nor}$ by the 
latter relation,
application of the trivialisation and translation from the left
and right.

The curvature $K$ inherits the  left- and right-invariance properties
from the canonical connection $\omega_{nor}$ and $K$ is also determined by 
its values
at the single point $(e,e)$ in $P\times N$. With respect to the 
trivialisation
$\iota_{\theta,g_{\frak{n}}}$, we find that
\begin{eqnarray*} E_i(\omega_{nor}(E_j))(e,e) &=& 
\frac{d}{dt}|_{t=0}\omega_{nor}(E_j)(\exp tE_i,e)=\frac{d}{dt}|_{t=0}
L^*_{\exp tE_i}\omega_{nor}(E_j)(e,e)\\
&=& \frac{d}{dt}|_{t=0}\omega_{nor}(E_j)(e,e)=0\end{eqnarray*}
for all $i,j=1,\ldots,n$. This shows for the curvature the identity
\[K(E_i,E_j)=-\omega_{nor}(e,e)([E_i,E_j]_{\frak{n}})
+[e_i+\gamma_{nor}(e_i),e_j+\gamma_{nor}(e_j)]_{\frak{g}}\ .\]
The curvature function $\kappa$ of the canonical Cartan connection 
$\omega_{nor}$ can then be expressed by (recall that the curvature is 
vertically trivial)
\[\kappa(e_i,e_j)=-(id+\gamma_{nor})\circ\rho_{\frak{n},\theta}(e_i,e_j)
+[e_i+\gamma_{nor}(e_i),e_j+\gamma_{nor}(e_j)]_{\frak{g}}\ .\]
Thereby, the $(-1)$-part of $\kappa$ is given through
\[\kappa_{-1}(e_i,e_j)=-\rho_{\frak{n},\theta}(e_i,e_j)+
[e_i,\gamma_0(e_j)]+
[\gamma_0(e_i),e_j]\ .\]
This expression vanishes, since $\omega_{nor}$ has no torsion. 
We see that the Lie bracket of $\frak{n}$ is
given on $\frak{m}_{-1}$ by 
\begin{eqnarray}\label{eq1}\rho_{\frak{n},\theta}(e_i,e_j)&=
&-\gamma_0(e_j)\cdot 
e_i+\gamma_0(e_i)\cdot e_j\ 
.\end{eqnarray}
The $0$-part of $\kappa$ is 
\[\kappa_0(e_i,e_j)=
-\gamma_0\circ\rho_{\frak{n},\theta}(e_i,e_j)+[e_i,\gamma_1(e_j)]+
[\gamma_1(e_i),e_j]+[\gamma_0(e_i),\gamma_0(e_j)]\ .\]
This part satisfies the trace-free condition
\begin{eqnarray} \label{eq2}\sum_{i=1}^n 
\gamma_0\circ\rho_{\frak{n},\theta}(e_i,a)(b)(e_i^*)=
\left\{\begin{array}{l}\quad\sum_{i=1}^n[e_i,\gamma_1(a)](b)(e_i^*)+
[\gamma_1(e_i),a](b)(e_i^*)\\[2mm] 
+\sum_{i=1}^n[\gamma_0(e_i),\gamma_0(a)](b)(e_i^*)
\end{array}\right\}\end{eqnarray} 
for all $a,b\in\frak{m}_{-1}$.
The $1$-part $\kappa_1$ of the curvature is
\[\kappa_1(e_i,e_j)=-\gamma_1\circ\rho_{\frak{n},\theta}(e_i,e_j)
+[\gamma_0(e_i),\gamma_1(e_j)]+[\gamma_1(e_i),\gamma_0(e_j)]\]
for all $i,j\in\{ 1,\ldots ,n\}$. 

The linear map $\gamma_{nor}:\frak{m}_{-1}\to\frak{p}$ is uniquely 
determined
be the normalisation conditions (\ref{eq1}) and (\ref{eq2}) with respect to
$\rho_{\frak{n},\theta}$ (which depends on the choice of the 
frame $\theta$).
(Otherwise, we would recover from another $\gamma$ which shares 
these properties
a further normal connection, which is not possible.) So $\gamma_{nor}$
depends only on the choice of $\theta$ which induces the Lie bracket
of $\frak{n}$ on $\frak{m}_{-1}$.
We can introduce the following formal notions.
\begin{DEF} Let $(G,P)$ be the flat homogeneous model (of conformal 
geometry) 
and 
\[\rho:\frak{m}_{-1}\times\frak{m}_{-1}\to\frak{m}_{-1}\]
a skew-symmetric map, which satisfies the Jacobi identity and  defines a 
Lie algebra 
bracket (of compact type) on the $(-1)$-part $\frak{m}_{-1}$ of the 
grading
of $\frak{g}$.
\begin{enumerate}\item 
We call a 
linear map
\[\gamma=\gamma_0+\gamma_1:\frak{m}_{-1}\to \frak{p}\]
(from the $(-1)$-part to the parabolic) 
a connection form on the model $(G,P)$. 
\item 
The curvature 
\[\kappa_{\gamma,\rho}
=\kappa_{-1}+\kappa_0+\kappa_1:\frak{m}_{-1}\times\frak{m}_{-1}
\to\frak{g}\] 
of the connection $\gamma$ with respect to the Lie bracket 
$\rho$ is defined as 
\[\kappa_{\gamma,\rho}(a,b)=-(id+\gamma)\circ\rho(a,b)
+[(id+\gamma)(a),(id+\gamma)(b)]_{\frak{g}}\] 
for $a,b\in\frak{m}_{-1}$. 
\item
The connection $\gamma$ is called torsion-free with respect to $\rho$
if $\kappa_{-1}=0$.
\item
The connection $\gamma$ is called normal with respect to $\rho$ if 
\[\kappa_{-1}=0\qquad\mbox{and}\qquad\mathrm{tr}\,\kappa_0=0\]
(cf. equation (\ref{eq2})).
\item There exists a unique normal connection with respect to 
the bracket $\rho$. We denote it by $\gamma_{\rho}:
\frak{m}_{-1}\to\frak{p}$ 
(or $\gamma_{nor}$ when the bracket is fixed on $\frak{m}_{-1}$) and 
call it the
canonical connection form of $\rho$ to the model $(G,P)$.
\end{enumerate}
\end{DEF} 

As we can see from formula (\ref{eq1}), the Lie bracket of $\frak{n}$
is determined on $\frak{m}_{-1}$ by the $0$-part $\gamma_0$ of the
normal connection $\gamma_{nor}$, since it has no torsion. In general, 
a skew-symmetric
map 
\[\gamma_0:\frak{m}_{-1}\to \frak{p}_0\] 
defines a Lie bracket on $\frak{m}_{-1}$ through
\[\rho_{\gamma_{0}}(a,b)=-\gamma_0(b)\cdot a+\gamma_0(a)\cdot 
b\qquad\mbox{for\ 
all\ }\ a,b\in \frak{m}_{-1}\]
if and only if the following sum of even permutations   
satisfies the relation (Jacobi identity)
\begin{eqnarray}\label{eq3}\sum_{\sigma(i,j,k)} 
[\gamma_0[e_i,\gamma_0(e_j)]+\gamma_0[\gamma_0(e_i),e_j]-
[\gamma_0(e_i),\gamma_0(e_j)],e_k]=0\end{eqnarray}
for all $i,j,k\in \{1,\ldots, n\}$.
The map $\gamma_0$ can then be extended in an arbitrary manner to a 
torsion-free connection form 
$\gamma$ with respect to $\rho_{\gamma_{0}}$ just by adding any linear 
$1$-part $\gamma_1$. In such a situation, the curvature function
to $\gamma$ with respect to $\rho_{\gamma_{0}}$ is given by
\begin{eqnarray*}\kappa_\gamma(e_i,e_j)&=&
-\gamma_0([e_i,\gamma_0(e_j)]+[\gamma_0(e_i),e_j])
+[e_i,\gamma_1(e_j)]+[\gamma_1(e_i),e_j]\\&&
+[\gamma_0(e_i),\gamma_0(e_j)]\\
&&-\gamma_1([e_i,\gamma_0(e_j)]+[\gamma_0(e_i),e_j])+
[\gamma_0(e_i),\gamma_1(e_j)]+[\gamma_1(e_i),\gamma_0(e_j)]\ 
.\end{eqnarray*}
Of course, not every torsion-free map 
$\gamma_0$ can be extended to the normal connection $\gamma_{nor}$ 
with respect to  $\rho_{\gamma_{0}}$. The condition on $\gamma_0$ for 
being normally extendible is given by the 
existence of a $\gamma_1$ such 
that
\begin{eqnarray}\label{eq4}
\sum_{i=1}^n([e_i,\gamma_1(a)]+[\gamma_1(e_i),a])(b)(e_i^*)=\!\!
\left\{\!\!\!\!\begin{array}{l}\quad\!\!\sum_{i=1}^n\!
\gamma_0([e_i,\gamma_0(a)]+[\gamma_0(e_i),a])(b)(e_i^*)\\[2mm]
-\sum_{i=1}^n[\gamma_0(e_i),\gamma_0(a)](b)(e_i^*)\end{array}\right. 
\!\!\!\end{eqnarray}
for all $a,b\in\frak{m}_{-1}$.

Finally, in this paragraph, we want to give explicit expressions
for the maps $\gamma_0$ and $\gamma_1$ of the normal connection
form $\gamma_{nor}$ and also for its curvature function $\kappa$  
when $\frak{n}$ is a semisimple Lie algebra (of 
compact 
type). The map $\gamma_0$ corresponds in this case to the Levi-Civita
connection of the bi-invariant metric $g_{\frak{n}}$ and is given
with respect to a reference frame $\theta$ by 
\[\gamma_0(e_i)=-\theta\circ\nabla_{\theta^{-1}(\cdot)} 
E_i=\frac{1}{2}\rho_{\frak{n},\theta}(e_i,\cdot)\]
for all $i=1,\ldots,n$.
Obviously, the so-defined map $\gamma_0$, considered as a matrix in 
$\frak{p}_0=\frak{co}(n)$ with respect
to the basis $(e_1,\ldots,e_n)$, satisfies (\ref{eq1}) (also 
(\ref{eq2})), i.e., 
$\gamma_0$
is torsion-free with respect to $\frak{n}$ (and $\theta$). 
Then we calculate for the traces on the left hand side in (\ref{eq4}):
\[\begin{array}{l}\sum_{i=1}^n 
\gamma_0([e_i,\gamma_0(a)]+[\gamma_0(e_i),a])(b)(e_i^*)
=\frac{1}{2}B_{\frak{n}}(a,b)\ ,\\[3.5mm]
\sum_{i=1}^n[\gamma_0(e_i),\gamma_0(a)](b)(e_i^*)=
\frac{1}{4}B_{\frak{n}}(a,b)\ .\end{array}\]
We set $\lambda=\frac{-1}{8(n-1)}$ and 
\[\gamma_1(a)=\lambda a^*\]
for all $a\in\frak{m}_{-1}$. Calculation of the right hand side in
(\ref{eq4}) gives 
\[\sum_{i=1}^n([e_i,\gamma_1(e_k)]+[\gamma_1(e_i),e_k])(e_l)(e_i^*)
=2\lambda(n-1)\delta_{kl}\]
for all $k,l\in\{1,\ldots,n\}$. Comparing both sides
of (\ref{eq4}) proves that the 
normal connection form
for $\frak{n}$
is determined to
\[\begin{array}{cccl}\gamma_{nor}:& \frak{m}_{-1}&\to&\frak{p}\ 
.\\[2mm]
&a&\mapsto& \frac{1}{2}\rho_{\frak{n},\theta}(a,\cdot)-
\frac{1}{8(n-1)}a^*\end{array}\]
The curvature functions $\kappa_{-1}$ and $\kappa_1$ vanish 
identically,
since there is no torsion and the Cotton-York tensor $C$ of 
$g_\frak{n}$ is trivial. The $0$-part $\kappa_0$ is given 
by the Weyl tensor $W$ of $g_\frak{n}$. It is
\[\kappa(a,b)=\kappa_0(a,b)=\theta^*W(\theta^{-1}(a),\theta^{-1}(b))\ 
.\]

\section{Canonical Cartan connection and holonomy}       
\label{ab5}                         

In this paragraph, we conduct a further discussion
of the connection form $\gamma_{nor}$. 
However, one should keep in mind that properties of $\gamma_{nor}$
depend on the choice of the trivialisation which comes from
the bi-invariant metric, and therefore $\gamma_{nor}$ should be 
considered as a 'metric object'. 
Nevertheless, we will derive 
an explicit formula for
its holonomy algebra, which is then a 'purely' conformal invariant. 

Let 
\[\gamma_{nor}=\gamma_0+\gamma_1:\frak{m}_{-1}\to\frak{p}\]
be the normal connection to some semisimple Lie algebra 
$\frak{n}$ (of compact type) with induced bracket
\[\rho_{\frak{n}}=[\ ,\ ]_{\frak{n}}\]
on $\frak{m}_{-1}$ (and with respect to 
$\sigma_{g_{\frak{n}}}$ and some $\theta$). 
The first observation here is the following. 
One can show that the image of 
$\frak{m}_{-1}$ 
under the $0$-part $\gamma_0$ of $\gamma_{nor}$
is the 
Lie subalgebra $Der(\frak{n})$ in $\frak{p}_0$ 
consisting of all derivations
of $\frak{n}\cong\frak{m}_{-1}$. Since $\frak{n}$ is semisimple, the 
Lie algebra
$\frak{n}$ itself is naturally isomorphic to $Der(\frak{n})$ by
\[\begin{array}{cccc}\mbox{ad}:&\frak{n}&\to& Der(\frak{n})\ .\\[2mm]
&a&\to& \mbox{ad}_a\end{array}\]
However, we will see that the map $\gamma_0:\frak{m}_{-1}\to\frak{p}_0$ 
is not a Lie algebra isomorphism with respect to $\rho_{\frak{n}}$
on $\frak{m}_{-1}$.
In detail, we can verify these statements as follows.
Remember that the normal connection form to $\frak{n}$ is given by
\[ \gamma_{nor}(a)=\frac{1}{2}\rho_{\frak{n}}(a,\cdot)
-\frac{1}{8(n-1)}a^*\ .\]
With this expression in mind, it is obvious that the 
Jacobi identity in $\frak{n}$ implies
\[\gamma_0(x)\cdot[a,b]_{\frak{n}}=[\gamma_0(x)\cdot a,b]+
[a,\gamma_0(x)\cdot b]\]
for all $x,a,b\in\frak{m}_{-1}$, where the 'dot' denotes matrix 
multiplication of $\frak{p}_0$ on $\frak{m}_{-1}$. 
This shows that 
$\gamma_0(x)$ 
is a derivation on $\frak{n}\cong(\frak{m}_{-1},\rho_{\frak{n}})$
for all $x\in\frak{m}_{-1}$. Moreover, the kernel of $\gamma_0$ is
trivial, since $\gamma_0(x)=0$ implies that $x$ is in the center of
$\rho_{\frak{n}}$, which itself is trivial for semisimple $\frak{n}$. We 
can
conclude that $\gamma_0$ is a vector space isomorphism onto 
$Der(\frak{n})$ and,
since 
$B_{\frak{n}}([a,b]_{\frak{n}},b)=0$ for $a,b\in\frak{n}$,         
every derivation sits in the semisimple part $\frak{so}(n)$ of
$\frak{p}_0$.
In fact, the following statement is true.
\begin{PR} 
Let $\frak{n}$ be a semisimple Lie algebra (of compact
type) and $\gamma_{nor}=\gamma_{\frak{n}}$ the corresponding connection
form on the model $(\frak{g},\frak{p})$ (coming from $(G,P)$).
Then 
the $0$-part $\gamma_0$
of $\gamma_{nor}$ is with respect to $\frak{n}$ the only torsion-free 
map, which sends $\frak{m}_{-1}$ to the derivations $Der(\frak{n})$
sitting in $\frak{p}_0$.\label{pr1} 
\end{PR}
{\bf Proof.} As we have seen, the map  $\gamma_0$ of $\gamma_{nor}$ admits
the stated properties. We have to 
show that any other linear map $Q$ apart from $\gamma_0$
sending $\frak{m}_{-1}$ to the derivations $Der(\frak{n})$ 
has torsion. This can be seen as follows. Since $Q$ is 
linear and $\gamma_0$ onto
$Der(\frak{n})$, we can write $\gamma_0-Q=\gamma_0\circ A$ for some 
homomorphism $A$ on $\frak{m}_{-1}$. The map $Q$ has 
no torsion if and only if 
\[-\gamma_0\circ A(e_j)e_i+\gamma_0\circ A(e_i)e_j=0\qquad
\mbox{for\  all\ }\ i,j\ .\]
Since $\gamma_0\circ A(a)\in\frak{so}(n)$ for all $a\in\frak{m}_{-1}$,
we find in the torsion-free case the relation
$\gamma_0\circ A(a)a=0$ for all $a\in\frak{m}_{-1}$. But this 
also implies \[\gamma_0\circ A(e_j)e_i+\gamma_0\circ A(e_i)e_j=0\qquad
\mbox{for\  all\ }\ i,j\in1,\ldots,n\ ,\]
which is only possible if $A=0$, i.e., $Q=\gamma_0$.\hfill$\Box$

Proposition \ref{pr1} is a characterisation of the 
Levi-Civita connection to the bi-invariant metric. The stated 
criteria replace conditions (\ref{eq3}) and (\ref{eq4}) for a map 
$\gamma_0$ to be the $0$-part of the normal connection $\gamma_{nor}$
for some semisimple $\frak{n}$.

The $0$-part $\kappa_0$ of the
curvature function to $\gamma_{nor}$ of $\frak{n}$ satisfies the 
relation
\[\kappa_0(a,b)-([a,\gamma_1(b)]+
[\gamma_1(a),b])=   
-\gamma_0\circ\rho_{\frak{n}}(a,b)
+[\gamma_0(a),\gamma_0(b)]\ .\]
The expression on the left hand side of this equation is the Riemannian
curvature tensor $R|_{\frak{n}}$ for the bi-invariant metric 
$g_{\frak{n}}$
in the conformal class $c_{\frak{n}}$. Obviously, it 
measures the deviation of
$\gamma_0$ from
being a Lie algebra isomorphism onto $Der(\frak{n})$ sitting
in $\frak{p}_0$. Since
the Riemannian curvature tensor of a semisimple Lie algebra is not 
zero, we can conclude that $\gamma_0$ is never a Lie algebra 
isomorphism. 
We state the following summary about the meaning 
of the curvature function 
$\kappa_{\gamma,\frak{n}}$ of 
an (arbitrary)
connection form $\gamma:\frak{m}_{-1}\to\frak{p}$ with respect to
some semisimple Lie algebra $\frak{n}$:
\begin{enumerate}\item
The $(-1)$-part $\kappa_{-1}$ of $\kappa_{\gamma,\frak{n}}$ measures 
the deviation
of the expression
\[-\gamma_0(b)a+\gamma_0(a)b\ ,\qquad a,b\in\frak{m}_{-1}\]
from being a defining map for the Lie bracket of 
$\frak{n}$ induced on $\frak{m}_{-1}$.
\item 
The $0$-part $\kappa_0$ is the deviation of the map $\gamma_0$ from
being a Lie algebra homomorphism onto $Der(\frak{n})$ up to an expression
depending on $\gamma_1$.
\item
The $1$-part $\kappa_1$ measures the deviation of the dual 
map \[\gamma_1^*:\frak{m}_{-1}\to\frak{m}_{-1}\]
to $\gamma_1$ from being a derivation on $\frak{m}_{-1}\cong\frak{n}$.
\item
In case that $\gamma=\gamma_{\frak{n}}$ is the normal connection
to $\frak{n}$ there is 
no torsion and $\gamma_1$ is a derivation on
$\frak{m}_{-1}$. The part $\gamma_0$ always differs from 
a homomorphism onto $Der(\frak{n})$. 
The connection $\gamma_{nor}$ is flat, i.e., $\kappa=0$, if and only if
$\frak{n}\cong\frak{su}(2)$.
\end{enumerate}

We come now to the second part of this paragraph concerning
the holonomy of the canonical Cartan connection. Recall that
the conformal automorphism
group of a compact Riemannian space is always an isometry group
with respect to some metric in the conformal class except in the case
when the compact space is the sphere $S^n$ with 
canonical 
conformally flat 
structure (M{\"o}bius sphere). In particular, the conformal automorphism
group of a compact semisimple Lie group with bi-invariant metric consists
entirely of the right and left translations with respect to the group 
multiplication
except for the case when $N=\SU(2)$. However, the normal Cartan connection 
on 
a conformal space gives rise to another
conformal invariant,
which can be thought of as an expression for conformal symmetry 
on a space (and yet seems not to be naturally related to a particular
metric in the conformal class),
namely the holonomy group $Hol(\omega_{nor})$ of the canonical connetion
$\omega_{nor}$.
To define this
holonomy group, we use the natural extension of the normal 
Cartan connection $\omega_{nor}$ on the parabolic frame bundle
$P(M)$ over a space $M$ 
with conformal structure $c$ to a usual principal fibre bundle connection 
$\tilde{\omega}_{nor}$ on the extended bundle
\[G(M)=P(M)\times_{P}G\]
with structure group $G=\SO(1,n+1)$. 
\begin{DEF} Let $(M,c)$ be a space with conformal structure. 
The holonomy 
group $Hol(\omega_{nor})$ of the canonical Cartan connection
$\omega_{nor}$
is defined to be the holonomy group of the naturally extended principal 
fibre 
bundle 
connection $\tilde{\omega}_{nor}$ on $G(M)$, whose group elements arise 
in the usual way by parallel translation of a fibre in $G(M)$ along closed 
curves on $M$.
The Lie algebra of $Hol(\omega_{nor})$ is denoted by 
$\frak{hol}(\omega_{nor})$. 
\end{DEF}
We remark that  there is a direct way to define the holonomy group of a 
Cartan 
connection without using the extended bundle. However, this 
approach results in general 
to the same group as in our definition
(cf. \cite{Sha97}). Note also that the holonomy group 
$Hol(\omega_{nor})$ is always a 
closed
subgroup of the M{\"o}buis group $G=\SO(1,n+1)$. We want to derive
here a formula, which can be used for explicit calculations
of the conformal holonomy on Lie groups with bi-invariant metric.
  
Let $N$ be a connected and compact semisimple Lie group. Let $\frak{n}$
denote its Lie algebra with bi-invariant Riemannian metric $g_\frak{n}$.
The canonical Cartan connection $\omega_{nor}$ induces in a natural
manner (with respect to the trivialisation coming from $g_{\frak{n}}$)
the map 
\[\gamma_{nor}:\frak{m}_{-1}\cong\frak{n}\to\frak{p}\ ,\]
which possesses all the informations of $\omega_{nor}$ and, in fact, can
be used to recover the canonical Cartan connection on $P(N)$.
We denote by
\[\Lambda(\frak{m}_{-1}):=span\{(id+\gamma_{nor})(a)|\
a\in\frak{m}_{-1}\}\subset \frak{g}\]
the images of the normal connection and by
\[\frak{q}:= span\{\kappa_{\frak{n}}(a,b)|\ 
a,b\in\frak{m}_{-1}\}\subset\frak{p}\]
the vector space of curvature values to the connection $\gamma_{nor}$.
There is a classical
formula for the holonomy algebra of an invariant connection 
on a homogeneous space with arbitrary structure group $G$ (cf. 
\cite{KN63}). 
We use this result to derive here easily the following formula for the 
conformal holonomy algebra of a bi-invariant metric.
\begin{THEO} (cf. \cite{KN63}) \label{TH1} Let $N$ be a connected and 
compact 
semisimple Lie group
with conformal structure $[g_\frak{n}]$. Then the holonomy algebra 
of the normal Cartan connection $\omega_{nor}$ on $(N,[g_\frak{n}])$
is given by the expression
\[\frak{hol}(\omega_{nor}):=
\frak{q}+[\Lambda(\frak{n}),\frak{q}]+[\Lambda(\frak{n}),
[\Lambda(\frak{n}),\frak{q}]]+\cdots\ ,\]
which is a subalgebra of $\frak{g}=\frak{so}(1,n+1)$.
The reduced holonomy group $Hol_o(\omega_{nor})$ is the connected 
subgroup 
of the identity component $SO_o(1,n+1)$ belonging to 
$\frak{hol}(\omega_{nor})$.
\end{THEO}

Usually, we do all our calculations with respect to the bi-invariant
metric $g_{\frak{n}}$. According to this,  we want to present the above 
formula for the conformal
holonomy algebra in a more 'suggestive' form. Let us denote
\[\frak{LC}:=span\{\gamma_0(a)|\ a\in\frak{m}_{-1}\}\ \subset\frak{p}_0\ 
.\]
This space is the span of the values in $\frak{so}(n)$ generated by the 
Levi-Civita connection of $\frak{g}_{n}$ and is isomorphic 
as Lie algebra to the derivations $Der(\frak{n})$. Then it is
\[\Lambda(\frak{m}_{-1})=
span\{\gamma_0^{-1}(l)+l+\gamma_1\circ\gamma_0^{-1}(l)|\ 
l\in\frak{LC}\}\ .\]
In short, we use the notation $\Lambda(\frak{m}_{-1})=\gamma_{nor}^{-1}
(\frak{LC})$. This space is isomorphic to $\frak{LC}$ as vector space,
but it is not a subalgebra in $\frak{g}$.
Moreover, we set
$\frak{W}:=\frak{q}$, which just shall remember to the fact
that $\frak{q}$ is generated from the values
of the Weyl tensor $W$ of $\frak{g}_{n}$. Then our formula for the
holonomy algebra takes the form
\[\frak{hol}(\omega_{nor})=\frak{W}+[\gamma_{nor}^{-1}(\frak{LC}),\frak{W}]
+[\gamma_{nor}^{-1}(\frak{LC}),[\gamma_{nor}^{-1}(\frak{LC}),\frak{W}]]+\cdots\ 
.\]
To compare, the holonomy algebra of the Levi-Civita connection
of the bi-invariant metric $g_{\frak{n}}$ is given by
\[\frak{hol}(g_{\frak{n}})=\frak{R}+[\frak{LC},\frak{R}]
+[\frak{LC},[\frak{LC},\frak{R}]]+\cdots\
,\]
where $\frak{R}$ denotes the space of images of the Riemannian 
curvature tensor to $g_\frak{n}$.

\section{Examples}           
\label{ab6}                  

We want to make some explicit use of the developed
calculus for conformal geometry on bi-invariant metrics.

{\bf Example 1.} 
Let $N=\SO(3)$ be the special orthogonal group in dimension
$3$, which is a $3$-dimensional compact and semisimple Lie group. 
Let $\frak{so}(3)$ denote its Lie algebra. We use 
for $\frak{so}(3)$ the 
standard basis $\{E_{ij}|\ 1\leq i<j\leq 3\}$, were the $E_{ij}$'s
are defined by matrix operations as
\[ E_{ij}:= e_i\cdot {}^te_j-e_j\cdot{}^te_i\]   
with respect to the standard basis $(e_1,e_2,e_3)$ of $\RR^3$.

The Lie algebra $\frak{so}(3)$ is isomorphic to $\frak{su}(2)$ and the 
universal covering of the group $\SO(3)$ is 
\[S^3=\Spin(3)\cong\SU(2)\ .\]
The bi-invariant metric on $\SO(3)$ is conformally
flat, since the Weyl tensor $W$ vanishes in dimension $3$ in general,
and $C$ vanishes for any bi-invariant metric. Of course, 
this is also clear from the fact that the bi-invariant metric
on the universal covering group $\SU(2)$ is the standard 
metric on $S^3$. For that reason, the calculations of conformal curvature 
and holonomy should produce trivial results here.
The connection form $\gamma_{nor}=\gamma_0+\gamma_1$ can be presented 
in the following form.

The chosen basis $\{E_{ij}\}$ in $\frak{so}(3)$ is orthogonal and 
$-B(E_{ij},E_{ij})=2$ for all its elements. We set the frame
\[\theta(\frac{1}{\sqrt{2}}E_{12})=e_1,\qquad
\theta(\frac{1}{\sqrt{2}}E_{13})=e_2,\qquad
\theta(\frac{1}{\sqrt{2}}E_{23}) =e_3\ .\]   
Then we find 
\[\begin{array}{l}\gamma_0(e_1)=-\frac{1}{\sqrt{2}}\nabla E_{12}=
\frac{1}{2\sqrt{2}}[E_{12},\cdot]=\frac{1}{2\sqrt{2}}E_{23}\ ,\\[4mm]
\gamma_0(e_2)=-\frac{1}{\sqrt{2}}\nabla E_{13}=
\frac{1}{2\sqrt{2}}[E_{13},\cdot]=-\frac{1}{2\sqrt{2}}E_{13}\ ,\\[4mm] 
\gamma_0(e_3)=-\frac{1}{\sqrt{2}}\nabla E_{23}=
\frac{1}{2\sqrt{2}}[E_{23},\cdot]=\frac{1}{2\sqrt{2}}E_{12}\ 
.\end{array}\] 
For $\gamma_1$ we have 
\[\gamma_1(e_i)=-\frac{1}{16}e_i^*\qquad i=1,2,3\  .\]
These identities express the normal Cartan connection on $\SO(3)$
in the trivialisation of the bi-invariant metric.
In fact, one can see from this $\gamma_{nor}$ that the curvature 
$\kappa_0$ vanishes identically. 
In particular, the conformal holonomy
group is trivial, which results from the formula in
Theorem \ref{TH1} with $\frak{q}=0$. 

To complete the discussion, we state that
the conformal automorphism
group of the covering space $\SU(2)$ is the M{\"o}bius group
$\SO(1,4)$. However, the conformal automorphism group
of $\SO(3)$ consists only of the isometry group 
\[\SO(4)/\ZZ_2=\SO(3)\times\SO(3)\]
of the bi-invariant metric.

{\bf Example 2.} 
We apply the conformal Cartan calculus now to the $6$-dimensional
compact and semisimple Lie group $\SO(4)$ with bi-invariant
metric induced by the Killing form $B$. The conformal automorphism
group of $\SO(4)$ consists entirely of the isometries, which
are the left and right translations, i.e. $\SO(4)\times \SO(4)$. 
The Lie algebra
$\frak{so}(4)$ is isomorphic to \[\frak{so}(3)\oplus\frak{so}(3)\]
and as its basis we use two copies 
of the basis $\{E_{ij}\}$ of $\frak{so}(3)$, namely
\[\{E_{ij},E_{kl}|\ 1\leq i<j\leq 3\ \ \mbox{and}\ \ 4\leq k<l\leq 6\}\ 
.\]
This basis is orthogonal with $-B(E_{ij},E_{ij})=2$
and it provides an embedding of $\frak{so}(4)$ into $\frak{so}(6)$.

The bi-invariant
metric $g_{\frak{so}(4)}$ induced by 
the Killing form on $\SO(4)$ is Einstein with positive 
scalar curvature. Obviously, it has non-constant sectional curvature
$S$ (cf. paragraph \ref{ab3}). Hence, it is not conformally flat. For that
reason, we expect in our calculation non-trivial curvature and holonomy  
for the conformal structure $c_{\frak{so}(4)}$ on $\SO(4)$. 

First, we calculate the normal connection form 
$\gamma_{nor}=\gamma_0+\gamma_1$. We use the frame
\[\theta(\frac{1}{\sqrt{2}}E_{12})=e_1,\qquad   
\theta(\frac{1}{\sqrt{2}}E_{13})=e_2,\qquad
\theta(\frac{1}{\sqrt{2}}E_{23}) =e_3\ ,\]
\[\theta(\frac{1}{\sqrt{2}}E_{45})=e_4,\qquad   
\theta(\frac{1}{\sqrt{2}}E_{46})=e_5,\qquad
\theta(\frac{1}{\sqrt{2}}E_{56}) =e_6\ .\]
From the calculations
for the case of $\frak{so}(3)$ we get 
\[\begin{array}{l}\gamma_0(e_1)=
\frac{1}{2\sqrt{2}}E_{23},\qquad 
\gamma_0(e_2)=-\frac{1}{2\sqrt{2}}E_{13},\qquad    
\gamma_0(e_3)=\frac{1}{2\sqrt{2}}E_{12},\\[4mm]
\gamma_0(e_4)=-\frac{1}{2\sqrt{2}}E_{56},\qquad
\gamma_0(e_5)=-\frac{1}{2\sqrt{2}}E_{46},\qquad
\gamma_0(e_6)=\frac{1}{2\sqrt{2}}E_{45}\ . \end{array}\]
The $1$-part $\gamma_1$ is given by
\[\gamma_1(e_i)=-\frac{1}{40}e_i^*\qquad i=1,\ldots, 6\ .\]

The $(-1)$-part $\kappa_{-1}$ of the curvature vanishes, since 
$\gamma_{nor}$
is torsion-free. The $1$-part $\kappa_1$ also vanishes,
since the Cotton-York tensor $C$ of $g_{\frak{so}(4)}$ is zero.
The $0$-part $\kappa_0$ consists of the Weyl tensor $W$.
As next we have to calculate the images of $\kappa_0$, which will give us 
the 
space $\frak{q}$.
It is 
\[W=R+\frac{1}{8(n-1)}g_\frak{n}*g_\frak{n}\ .\]
For the Kulkarni-Nomizu product in this sum it is easy to see
that
\[\theta^*(B*B)(\theta^{-1}(e_i),\theta^{-1}(e_j))
=\frac{1}{20}E_{ij}\qquad \mbox{for\ all\ 
}i,j=1,\ldots,6\ .\]
For the Riemannian curvature tensor we find 
\[\theta^*R_\frak{n}(\theta^{-1}(e_i),\theta^{-1}(e_j))
=-\frac{1}{8}E_{ij}\]
for all $i,j\in\{1,\ldots 3\}$ and $i,j\in\{4,\ldots 6\}$.
The remaining curvature expressions for $R_{\frak{n}}$ are zero.
This shows that the span of the $0$-part $\kappa_0$ of
the Cartan curvature is equal to $\frak{so}(6)$, which is the semisimple
part of $\frak{p}_0$ in the M{\"o}bius algebra $\frak{so}(1,7)$.
Hence, we get for the span of the curvature values
\[\frak{q}=\frak{so}(6)\ \subset\frak{p}_0\ .\]
Obviously, the span of the normal connection $id+\gamma_{nor}$
is given by
\[\Lambda(\frak{m}_{-1})=\{e_i+\gamma_0(e_i)-\frac{1}{40}e_i^*|\ 
i=1,\ldots, 6\}\ .\]
Then it is a straightforward calculation to see that the 
space $[\Lambda(\frak{m}_{-1}),\frak{q}]$ 
of commutators is 
equal to
\[span\{e_i-\frac{1}{40}e_i^*|\ i=1,\ldots,6\}\oplus
\frak{so}(6)\ .\]
We denote 
\[\frak{l}:=span\{x+\gamma_1(x)|\ x\in\frak{m}_{-1}\}\ .\]
The space $\frak{l}$ is stable under the action of $\frak{so}(6)$
sitting in $\frak{p}_0$. This shows that all the spaces 
\[[\Lambda(\frak{m}_{-1}),\cdots,[\Lambda(\frak{m}_{-1}),\frak{q}]\cdots]\]
of commutators
are equal to 
$\frak{l}\oplus\frak{q}$, which is seen to be  isomorphic to the Lie 
algebra 
$\frak{so}(7)$ 
embedded into $\frak{so}(1,7)$.
We conclude for the holonomy algebra that 
\[\frak{hol}(\omega_{nor})=\frak{l}\oplus\frak{q}\cong\frak{so}(7)\ .\]
The holonomy group of $\omega_{nor}$ on 
$\SO(4)$ is then given by
\[Hol(\omega_{nor})=\SO(7) .\]

This result has the following interpretation.
The bi-invariant metric $g_{\frak{n}}$ in the conformal class
$c_{\frak{n}}$ on $N=\SO(4)$ is Einstein. It is well known in general 
that the 
conformal Einstein condition 
implies that the holonomy group $Hol(\omega_{nor})$ of the normal Cartan
connection stabilises a standard 'tractor', that is a vector in the 
standard 
representation $\RR^{1,7}$ of the M{\"o}bius group $\SO(1,7)$ (cf.
e.g. \cite{Lei04a}).
In case that the scalar curvature of the Einstein metric is positive,
this vector is timelike, which explains that for our case the holonomy
group $Hol(\omega_{nor})$ of $\SO(4)$ is automatically reduced
to $\SO(7)$. Our calculation then shows that the holonomy
is not further reduced and we can say that $\SO(4)$ 
has generic conformal holonomy up to the fact that it is (conformally)
Einstein. In particular, we can read off from the holonomy result
that $\SO(4)$ does not admit any conformal Killing spinors nor
normal conformal Killing forms. Also, there is no
further Einstein metric in the conformal class $c_{\frak{n}}$ beside
the bi-invariant metric $g_{\frak{n}}$ (cf. \cite{Lei04a}).


\end{sloppypar}

\begin{thebibliography}{1111111}

\bibitem[Ber55]{Ber55} M. Berger. {\it Sur les groupes 
d'holonomie homogene des varietes a connexion affine et des varietes 
riemanniennes}, Bull. Soc. Math. France 83(1955), p. 279-330.
\bibitem[Lei04b]{Lei04b} F. Leitner,
{\it Canonical holonomy of a homogeneous parabolic geometry},
to appear, 2004.
\bibitem[O'N83]{O'N83} B. O'Neill. {\it 
Semi-Riemannian geometry}, Pure and Applied Mathematics. Academic Press, 
1983.
\bibitem[Sha97]{Sha97} R.W. Sharpe. {\it Differential Geometry},
Graduate Texts in Mathematics 166. Springer-Verlag New York, 1997
\bibitem[KN63]{KN63} S. Kobayashi, K. Nomizu. {\it 
Foundations of differential 
geometry}, John Wiley \& Sons, vol.1, New York, 1963.
\bibitem[Lei04a]{Lei04a} F. Leitner. {\it Normal conformal Killing 
forms}.
SFB288-Preprint No. 605, Berlin 2004.  
\bibitem[Kob72]{Kob72} S. Kobayashi. {\it Transformation Groups in
Differential Geometry}. Springer-Verlag Berlin Heidelberg, 1972.
\bibitem[CSS97]{CSS97}A. $\check{\rm{C}}$ap, J. Slovak,
V. Sou$\check{\rm{c}}$ek. {\it Invariant Operators on Manifolds with
Almost
Hermitian Symmetric Structures I \& II}. Acta. Math. Univ. Comen.,
New Ser.
66, No. 1, p. 33-69 \& No. 2, p. 203-220(1997).
\end{thebibliography}
\end{document}